\documentclass{article}

\usepackage{epsfig,amsmath,amsfonts}
\usepackage{amssymb}

\def\be{\begin{equation}}
\def\ee{\end{equation}}

\begin{document}
\title{The accurate numerical solution of the Schr\"odinger equation with an explicitly time-dependent Hamiltonian}
\author{V. Ledoux\footnote{Postdoctoral Fellow of the Fund for Scientific Research
- Flanders (Belgium) (F.W.O.-Vlaanderen)}, M. Van Daele}

%
\maketitle

\begin{abstract}
We show how the highly accurate and efficient Constant Perturbation (CP) technique for steady-state Schr\"odinger problems can be used in the solution of time-dependent Schr\"odinger problems with explicitly time-dependent Hamiltonians, following a technique suggested by Ixaru in \cite{Ixaru2010}.
By introducing a sectorwise spatial discretization using bases of accurately CP-computed eigenfunctions of carefully-chosen stationary problems, we deal with the possible highy oscillatory behaviour of the wave function while keeping the dimension of the resulting ODE system low. Also for the time-integration of the ODE system a very effective CP-based approach can be used.
\end{abstract}
%
%


\section{Introduction}
The constant perturbation methods form a class of numerical methods which were originally devised for the solution of the time-independent one-dimensional Schr\"odinger equation.  The first ideas were described in the seventies (see \cite{Ixaru1984} for references) but the CP methods still prove to be very efficient for not only the one-dimensional time-independent Schr\"odinger equation \cite{Ixaru1997} but also for Sturm-Liouville problems \cite{SLCPM12,Ledoux2005a},
coupled channel Schr\"odinger equations \cite{Ixaru2002,Ledoux2012} and two-dimensional Schr\"odinger problems \cite{Ixaru2010}. As noted in \cite{Degani}, the CP approach is very closely related to the application of a modified Neumann method. Modified Neumann methods, as well as CP methods, are very effective for highly oscillatory ordinary differential equations, see \cite{Neumann}.
In the current paper we extend the techniques to the solution of time-dependent Schr\"odinger problems involving explicitly time-dependent Hamiltonians. 

We consider the time-dependent Schr\"odinger equation  (in units where $\hbar=1$):
\begin{equation}
\imath\frac{\partial \psi(x,t)}{\partial t}={\hat H}(x,t)\psi(x,t),\quad x\in \mathbb{R},\; t>0
\label{eq0}
\end{equation}
with the one-dimensional time-dependent Hamiltonian
\[{\hat H}(x,t)=-\frac{1}{2\mu}\frac{\partial^2}{\partial x^2}+V(x,t),\]
and 
\[\psi(x,0)=\psi_0(x)\]
as the wave function at initial time $t=0$.
$\mu$ is the reduced mass. 

Standard techniques found in literate reduce \eqref{eq0} to the solution of a linear ordinary differential equation in $t$
of the form
\begin{equation}\imath\frac{d}{dt}{Z}(t)={H}(t){Z}(t),\label{eq1}\end{equation}
where ${Z}$ is a column vector with $N$ complex components and ${H}$ is an $N\times N$ hermitian matrix associated to the Hamiltonian.
Equation \eqref{eq1} arises when the spatial variable $x$ is discretized. Then the entries of ${Z}(t)$ are values of the wave function $\psi(x,t)$ at the nodes of the spatial grid. Equation \eqref{eq1} can also result from a spectral decomposition in which the solution  $\psi(x,t)$  is expanded over $N$ eigenfunctions of a time-independent Hamiltonian, e.g.\ a harmonic oscillator, and the ${Z}(t)$ represent the coefficients in this basis set expansion. 

Instead of taking a fixed set of basis eigenfunctions of a time-independent Hamiltonian over the entire $t$-range, we will use similar ideas as those proposed in \cite{Ixaru2010} to tackle 2D stationary equations in an efficient way. In \cite{Ixaru2010} Ixaru introduced a sectorwise expansion of the 2D solution over eigenfunctions from conveniently tuned 1D Hamiltonians.
In a similar way, we divide here the $t$-domain in sectors and use in each sector a different time-independent Hamiltonian to form the set of eigenfunctions. By choosing in each sector a time-independent Hamiltonian which takes into account the form of the potential $V(x,t)$ over this sector, the number $N$ of time-independent stationary solutions needed for the expansion can be kept low. The efficient solution of stationary Schr\"odinger problems forms a crucial part in the proposed spatial discretization. To compute the eigenvalues and the corresponding eigenfunctions, we make use of the highly accurate and efficient CP methods and codes which are available for solving time-independent Schr\"odinger problems.
Also for the propagation of the time-dependent solution $\psi$ over a sector, we will use CP-based techniques. 
By using a CP approach we are able to tackle the oscillations which can be present both in time and space without the necessity of small steps.

\section{The procedure}
We take the time-dependent equation \eqref{eq0} and assume that the $x$-range can be restricted to $[x_{min},x_{max}]$ and that Dirichlet boundary conditions can be imposed in the endpoints $x_{min}$ and $x_{max}$.

\begin{figure}
	\centering
		\includegraphics[width=5cm]{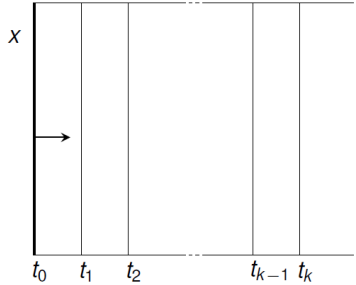} 
		\caption{The rectangular subdivision of the integration domain.}\label{fig1}
\end{figure}

\begin{figure}
	\centering
		\includegraphics[width=6cm]{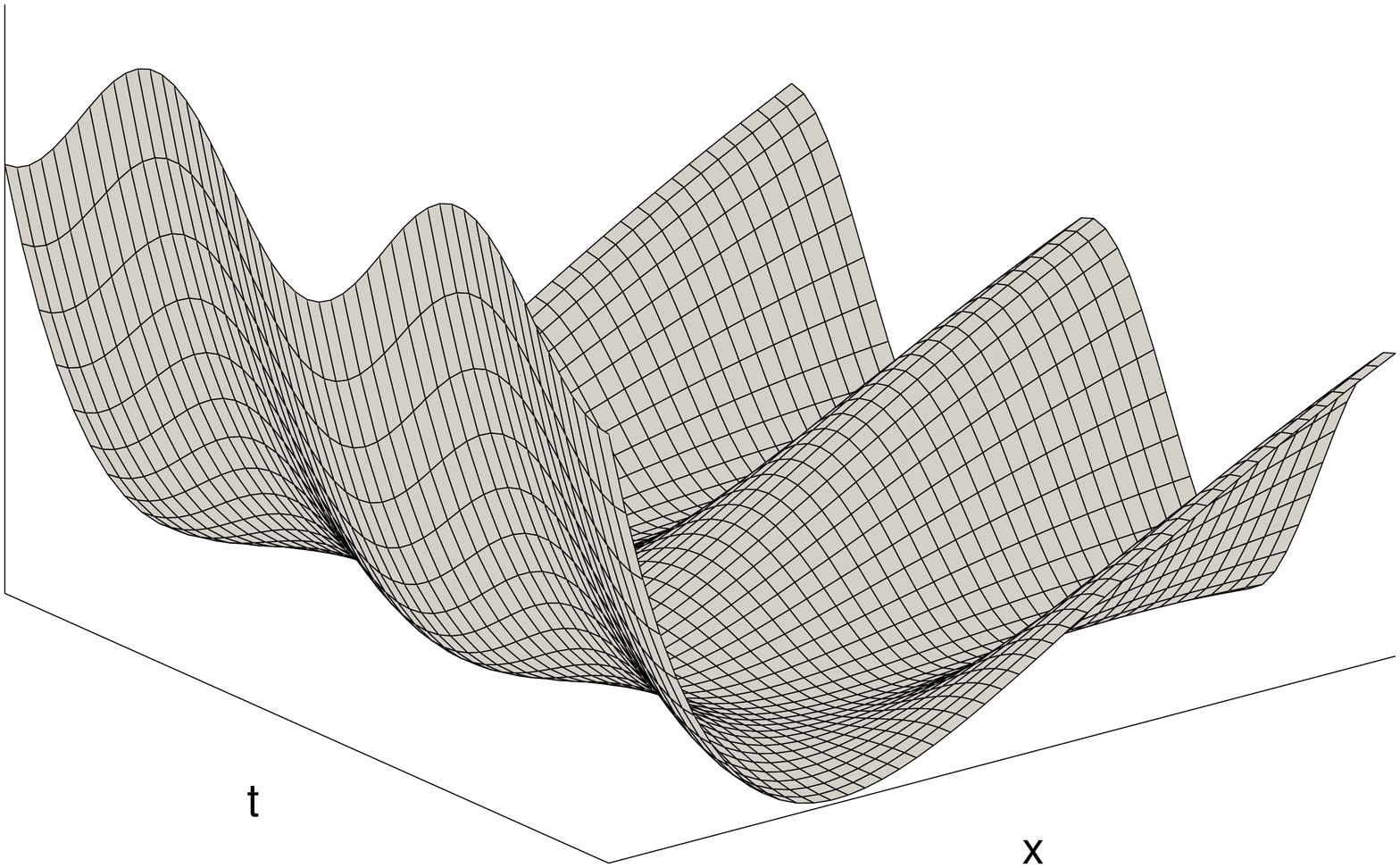} 
		\includegraphics[width=6cm]{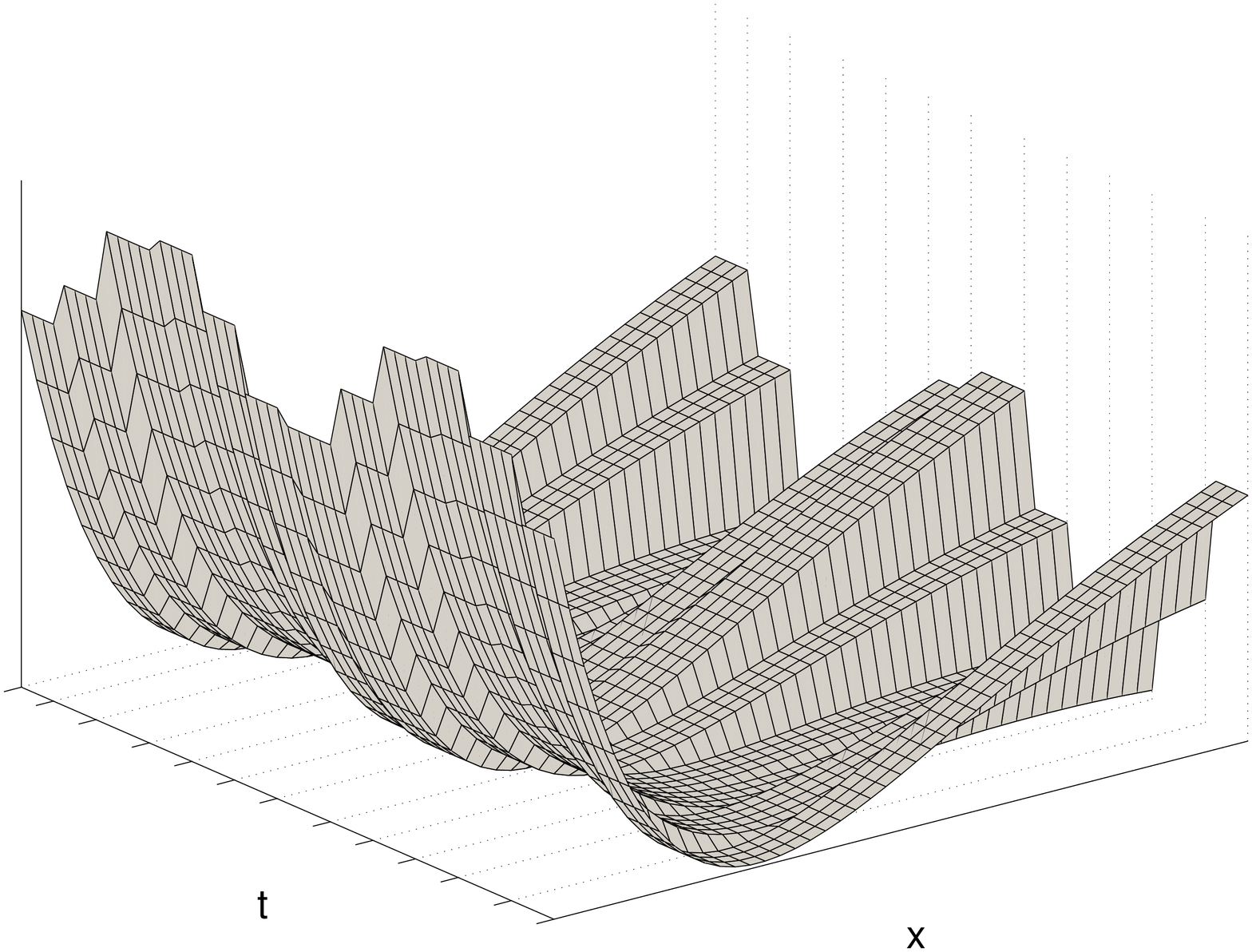} 
		\caption{The shape of a (Walker-Preston) time-dependent potential $V(x,t)$ and its sectorwise approximation ${\bar V}(x,t)$.}\label{fig2}
\end{figure}

First the $(x,t)$-domain is divided in a succession of rectangular sectors as illustrated in Figure \ref{fig1}. The mesh points on the $t$-axis need not to be equidistant. On each rectangular subdomain $[x_{min},x_{max}]\times[t_{k-1},t_k]$ the original time-dependent potential function $V(x,t)$ is approximated by the time-independent function ${\bar V}^{[k]}(x)=V(x,(t_{k-1}+t_k)/2)$. When doing this approximation over each sector, we obtain an approximation of the original $V(x,t)$ with a staircase shape in the $t$-direction, see Figure \ref{fig2}. 

Let us focus on the propagation of the solution over one sector. Assume the solution $\psi(x,t_{k-1})$ is known, either through the initial condition when $k=1$ or through the computation of the solution over the previous sector. First we solve the time-independent Schr\"odinger eigenvalue problem with potential function ${\bar V}^{[k]}(x)$:
\begin{align}\frac{1}{2\mu}\frac{d^2}{d x^2}y^{[k]}_n=\left({\bar V}^{[k]}(x)-E^{[k]}_n\right)y^{[k]}_n,\quad
y^{(k)}_n(x_{min})=y^{(k)}_n(x_{max})=0,\label{SE}\end{align}
where the upper label $[k]$ refers to the specific sector over which we are propagating. The eigenfunctions are normalized such that \be\left\langle y^{[k]}_n,y^{[k]}_m\right\rangle=\int_{x_{min}}^{x_{max}}y^{[k]}_n(x)y^{[k]}_m(x){\rm d}x=\delta_{nm}, \quad n,m=1,2,3,\dots .\label{eq2}\ee
The set of eigenfunctions $y^{[k]}_n$ form an orthonormal basis. We write the desired $\psi(x,t)$ over the $k$th sector as an expansion over this set, with $t$-dependent coefficients:
\be\psi^{[k]}(x,t)=\sum_{m=1}^{\infty}c_m^{[k]}(t)y^{[k]}_m(x).\label{eq3}\ee
Expansion \eqref{eq3} is then introduced into the time-dependent Schr\"odinger equation \eqref{eq0}. After multiplying this by $y^{[k]}_n(x)$, integrating over $x$ and using Eq. \eqref{eq2}, we obtain
\be \imath \frac{d}{dt}c_n^{[k]}(t)=\sum_{m=1}^{\infty}\left[\Delta V^{[k]}_{nm}(t)+E^{[k]}_n\delta_{nm}\right]c_m^{[k]}(t),\quad n=1,2,3,\dots\label{Csys}\ee
where
\be\Delta V^{[k]}_{nm}(t)=\int_{x_{min}}^{x_{max}} y^{[k]}_n(x) \left[V(x,t)-{\bar V}^{[k]}(x)\right]y^{[k]}_m(x){\rm d}x.\label{eq5}\ee
The propagation over sector $k$ is then reduced to the solution of the linear system \eqref{Csys} over the interval $[t_{k-1},t_k]$ with initial values $c_m^{[k]}(t_{k-1})$ which are derived from the data obtained on the previous sector. Continuity is imposed in $t_{k-1}$, that is
\[\sum_{m=1}^{\infty}c_m^{[k-1]}(t_{k-1})y^{[k-1]}_m(x)=\sum_{m=1}^{\infty}c_m^{[k]}(t_{k-1})y^{[k]}_m(x).\]
After multiplying this equation with $y^{[k]}_n(x)$ and integrating over $x$, we obtain the initial conditions:
\be c_n^{[k]}(t_{k-1})=\sum_{m=1}^{\infty} s_{nm}^{[k]}c_m^{[k-1]}(t_{k-1})\label{icc}\ee
where $s_{nm}^{[k]}$ measures the overlap between the eigenfunctions in sector $k$ and the ones in sector $k-1$:
\be s_{nm}^{[k]}=\int_{x_{min}}^{x_{max}}y^{[k]}_n(x)y^{[k-1]}_m(x){\rm d}x.\label{int2}\ee

In practice, some upper value $N$ must be imposed for the number of eigenfunctions included in the expansion \eqref{eq3}, and thus for the indices $n$ and $m$. This means that over each sector, $N(N+1)/2$ integrals need to be computed to obtain $\Delta^{[k]}V_{nm}(t)=\Delta^{[k]}V_{mn}(t)$ for $n,m=1,2,\dots,N$, and $N^2$ integrals to obtain each $s_{nm}^{[k]}$.
Also the system \eqref{Csys} can be written in matrix form as
\be{{C}^{[k]}}'(t)=A^{[k]}(t){C}^{[k]}(t),\quad t\in[t_{k-1},t_k],\label{system}\ee
where $C^{[k]}$ is a column vector with $N$ components $c_n^{[k]}$ and $A^{[k]}(t)$ is an $N$ by $N$ symmetric matrix with elements $A_{nm}^{[k]}(t)=-\imath \left[\Delta V^{[k]}_{nm}(t)+E^{[k]}_{n}\delta_{nm}\right]$. The initial condition for the system \eqref{system} is given by Eq.\ \eqref{icc} or in matrix form ${{C}^{[k]}}(t_{k-1})=S^{[k]}C^{[k-1]}(t_{k-1})$.

\section{Numerical methods}
The procedure described in the previous section reduces the solution of the time-dependent problem to the solution of time-independent Schr\"odinger problems \eqref{SE}, linear ODE systems of the form \eqref{system}, and the computation of the integrals \eqref{eq5} and \eqref{int2}.

\subsection{The time-independent 1D Schr\"odinger eigenvalue problem}
It is crucial to have a computationially efficient approximation technique for the solution of the time-independent Schr\"odinger problems \eqref{SE}, and to obtain a uniform high accuracy approximation to the eigenvalues and eigenfunctions. It is well known that the oscillatory behavior of the solutions of time-independent Schr\"odinger problems \eqref{SE}, force naive integrators to take increasingly smaller steps. As shown in \cite{Ixaru1997}, highly accurate and efficient techniques are available in the form of a shooting method employing Constant based Perturbation (CP) methods for the solution of the initial value problem. CP methods are based on ``approximation of the potential''. In particular, constant reference potential approximations are used and some perturbation corrections are computed from the difference between the constant reference potential and the original one. 

In our numerical experiments, we used a CP method of order 12 as it was implemented in the Matlab package {\sc Matslise} \cite{Ledoux2005a} (and earlier in Fortran form in \cite{SLCPM12}) to compute accurate approximations of the eigenvalues and evaluations of the eigenfunctions and their first order derivative in a fixed set of meshpoints. The mesh is formed by subdividing the $x$-range in a number of equidistant steps and taking the 4 Lobatto nodes over each such step as the meshpoints. Choosing the meshpoints in this way, is convenient for the fast and accurate approximation of the integrals \eqref{eq5} and \eqref{int2} over each sector (see further).

\subsection{Linear system of first order ODEs}
For Schr\"odinger equations with an explicitly time-dependent Hamiltonian, a difficulty when constructing accurate and efficient propagators is that fast oscillations in the time-dependent fields apparently necessitate small time steps. 
We propose here an effective discretization method for the solution of the linear system \eqref{system}, which takes into account the oscillatory nature of the solution. We again use a CP approach, which is strongly related to the application of a modified Neumann scheme, see \cite{Degani}. 
An alternative is to employ a modified Magnus approach, however the computation of the matrix exponential of the Magnus series terms is exceedingly expensive, especially when one wants to employ higher-order Magnus approximations.

First step is to approximate the matrix $A^{[k]}$ in equation \eqref{system} by a constant (time-independent) matrix ${\bar A}^{[k]}$. For the construction of a second order algorithm, we can use ${\bar A}^{[k]}=A^{[k]}((t_{k-1}+t_k)/2)$ as a constant approximation. Note that in this case ${\bar A}^{[k]}$ is a diagonal matrix with ${\bar A}^{[k]}_{nm}=-\imath E_n^{[k]}\delta_{nm}$ and that the solution can be easily propagated by 
\[C^{[k]}(t)=\exp[(t-t_{k-1}){\bar A^{[k]}}]C^{[k]}(t_{k-1}), \quad t \in[t_{k-1},t_k].\]
When constructing higher order CP methods, we approximate $A^{[k]}$ by a higher order polynomial approximation, i.e.\ by a truncated series over the shifted Legendre polynomials:
\[A^{[k]}(t_{k-1}+\delta)\approx {\hat A}(\delta)=\sum_{m=0}^M A_m h^m P^*_m(\delta/h),\quad h=t_k-t_{k-1},\; \delta \in [0,h],\]
where the matrix weights $A_m$ are calculated by quadrature (see \cite{Ixaru1997,LedouxNeumann}). The symmetric matrix ${\bar A}^{[k]}=A_0$ is then diagonalized and let $D$ be the orthogonal diagonalization matrix, i.e.\ $A_0=DA_0^{D}D^T$. Our propagation algorithm for the solution takes then the following form
\[C^{[k]}(t)=DT^D(t-t_{k-1})D^TC^{[k]}(t_{k-1})\]
with
\[T^D(\delta)=\exp(\delta{A_0^D}) + P_1(\delta)+P_2(\delta)+...\]
where $P_1,P_2,\dots$ are perturbation corrections derived from the perturbation ${\hat A}^D-A_0^D$. In fact, by using perturbation theory (as in \cite{Ixaru1984}) one can show that the perturbations $P_q, (q=1,2,\dots)$ satisfy 
\be P_q'(\delta)=A_0^D P_q(\delta) + [{\hat A}^D(\delta)-A_0^D]P_{q-1}(\delta),\quad P_q(0)=0,\label{rec}\ee
where $P_0(\delta)=\exp(\delta{A_0^D})$. 
This means that $P_1$ has the following form
\[P_1(\delta)=\exp(\delta{A_0^D})N^1(\delta),\]
where \be N^1=\int_0^\delta \exp(-\delta_1{A_0^D})[{\hat A}^D(\delta_1)-A_0^D]\exp(\delta_1{A_0^D})\,{\rm d}\delta_1,\ee corresponds exactly to the first term in a modified Neumann integral series. Also the expressions for the higher order perturbations are equal to the corresponding term in the modified Neumann series. Since ${\hat A}^D(\delta)$ has a polynomial form, the integrals in the Neumann series terms can be computed analytically. For a fourth order algorithm for instance, it is sufficient to take $M=1$ and one perturbation correction $P_1$ where $N^1(h)$ is given by
\[N^1_{ij}(h)=\frac{(h\Delta_{ji}+2)+(h\Delta_{ji}-2)\exp(h\Delta_{ji})}{\Delta_{ji}^2} (A^D_1)_{ij}.\]
with $\Delta_{ji}=({A_0^D})_{jj}-({A_0^D})_{ii}$ and $A^D_1=D^TA_1D$.

In our implementations we used this simple fourth order scheme and applied it once over the full sector-length $[t_k,t_{k+1}]$ or (if needed to ensure accuracy) on a subdivision of the sector $[t_k,t_{k+1}]$. 

It is often important to discretize Schr\"odinger equations while preserving unitarity so that the conservation of the norm is ensured. The second order scheme always preserves the norm and non-unitarity appears only through the correction terms which are added. The loss of norm conservation is consequently likely to be small, moreover including correction terms improves the accuracy of the method, i.e.\ it leads to an approximate solution which is closer to the exact one. We illustrate this in the numerical experiments section.

\subsection{Oscillatory integrals}
Also for the quadrature of the integrals involved in the procedure it is important to select the appropriate rules. The integrands of \eqref{eq5} and \eqref{int2} contain a product of two wavefunctions of a time-independent Schr\"odinger equation. These wavefunctions have an oscillatory character in that part of the integration domain where the corresponding energy eigenvalue is larger than the potential.
An alternative to the standard quadrature methods is provided by contemporary methods for highly oscillatory quadrature, an
area that has undergone significant developments over the last ten years. Many of these methods were, however, specially devised for integrals of the form $\int_a^b f(x)e^{\imath \omega g(x)} {\rm d}x$ where $\omega\in \mathbb{R}, \omega \gg 1$ and the nonoscillatory functions $f$ and $g$ are known \cite{Huybrechs,IN,Olver2006,HO}. Eqs.\ \eqref{eq5} and \eqref{int2} do not fit in this $\int_a^b f(x)e^{\imath \omega g(x)} {\rm d}x$ framework. Moreover, it may be more convenient for the present application to use the same set of quadrature nodes for each integral instead of frequency-dependent ones, in order to minimize the number of CP eigenfunction evaluations.  However, we do have some information available which we can use in the formulation of a specially tuned quadrature rule: we have some knowledge about the frequency of each wavefunction, moreover the CP approach does not only provide us with the values of the wavefunctions at the meshpoints but also those of their first derivative. We will use this information in the construction of some adapted quadrature rules based on exponential fitting (EF). These EF rules are not only very effective in the oscillatory part of the integration domain but also in the region where an eigenfunction displays an exponential behavior. Details about EF quadrature can be found in \cite{Ixaru2004}.

We first consider integral \eqref{int2}. The integrand can be generically written as the product $I(x)=u(x)z(x)$, where $u$ is an eigenfunction of the Schr\"odinger problem $u''(x)=(Q_u(x)-\lambda_u)u(x)$ and $z$ an eigenfunction of a different Schr\"odinger problem $z''(x)=(Q_z(x)-\lambda_z)z(x)$. As described in \cite{Ixaru2004}, p.\ 47, the eigenfunction $u(x)$ corresponding to an eigenvalue $\lambda_u=2\mu E_u$, can be written over the meshinterval $[x_i,x_{i+1}]$ as
\[u(x)=f_1(x)\exp(\sqrt{{\bar Q}_u-\lambda_u} x)+f_2(x)\exp(-\sqrt{{\bar Q}_u-\lambda_u} x)\]
where ${\bar Q}_u$ is a constant approximation of the potential $Q_u(x)=2\mu {\bar V}^{[k]}(x)$ on the current meshinterval. Note that we have such a constant approximation of the potential available in the CP algorithm. Similarly,
\[z(x)=g_1(x)\exp(\sqrt{{\bar Q}_z-\lambda_z} x)+g_2(x)\exp(-\sqrt{{\bar Q}_z-\lambda_z} x).\]
Based on this structure for both $u$ and $z$, we choose to use a 4-point Lobatto-EF algorithm of the form
\be \int_{x_i}^{x_{i+1}} I(x){\rm d}x=\int_{X-h}^{X+h} I(x){\rm d}x\approx h\sum_{n=1}^{4}a_n^{(0)}I(X+x_nh)+h^2\sum_{n=1}^{4}a_n^{(1)}I'(X+x_nh)\label{Efrule}\ee
which is exact for the functions
\[\exp(\pm \mu_1 x),\exp(\pm \mu_2 x),x\exp(\pm \mu_1 x),x\exp(\pm \mu_2 x),\]
where $\mu_1=\sqrt{{\bar Q}_y-\lambda_y}+\sqrt{{\bar Q}_z-\lambda_z}$ and $\mu_2=\sqrt{{\bar Q}_y-\lambda_y}-\sqrt{{\bar Q}_z-\lambda_z}$. See \cite{Ixaru2004} for the construction of such an EF rule. 

The integrand of integral \eqref{eq5} has a bit different form which we generically write as $f(x)u(x)z(x)$. $u$ and $z$ are again eigenfunctions for which we have the evaluations of the first derivative available. In many cases the form of the potential function $V(x,t)$ is explicitly known and the first derivative w.r.t.\ $x$ of the function $f(x)=V(x,t)-{\bar V}^{[k]}(x)$ can be evaluated. In this case, the EF method \eqref{Efrule} is used. Otherwise, we propose an EF scheme which does not use derivative information and is exact for $\exp(\pm \mu_1 x),\exp(\pm \mu_2 x)$. 

\section{Numerical results}

\begin{figure}
	\centering
	\includegraphics[width=6cm]{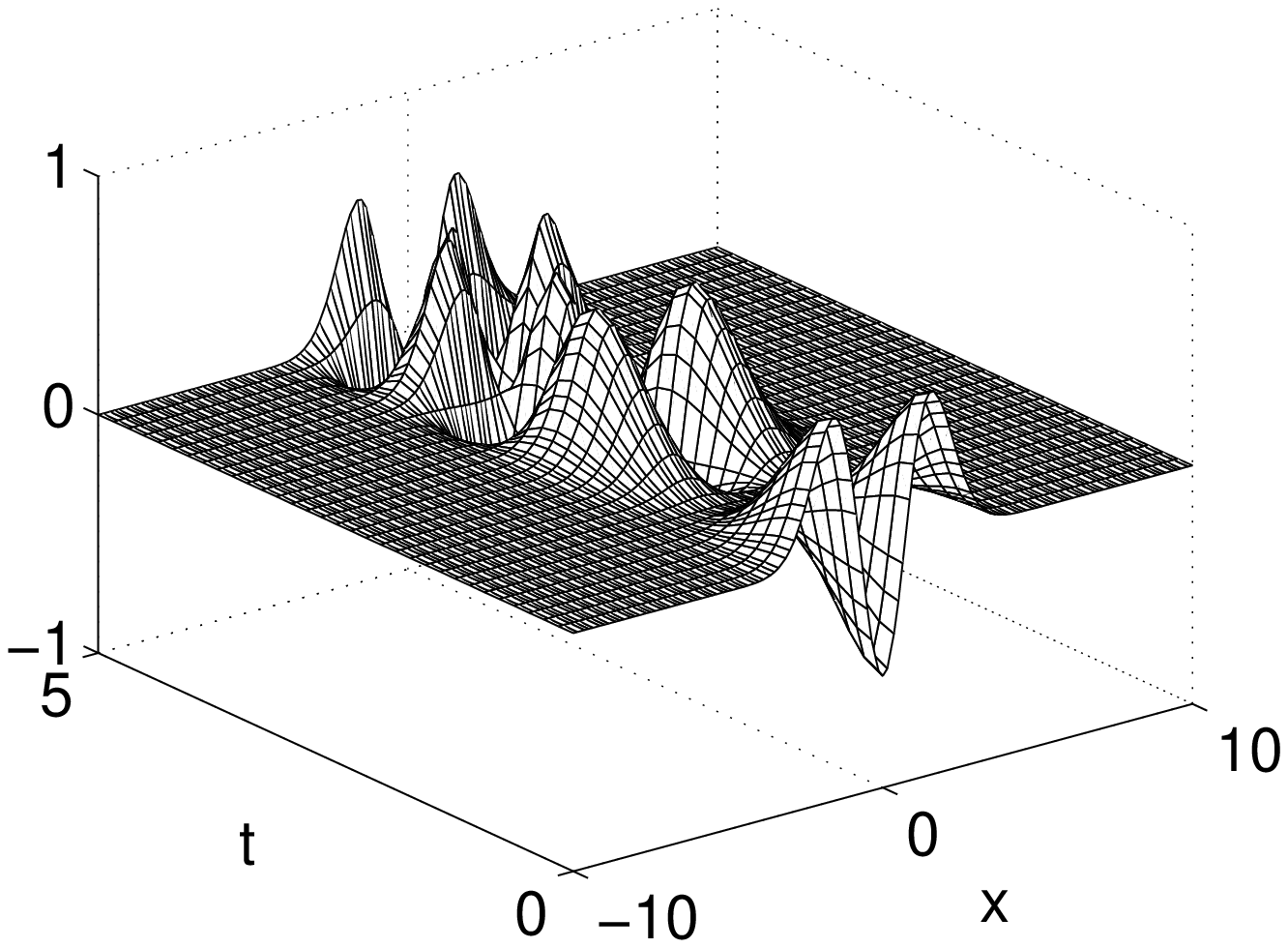}
	\includegraphics[width=6cm]{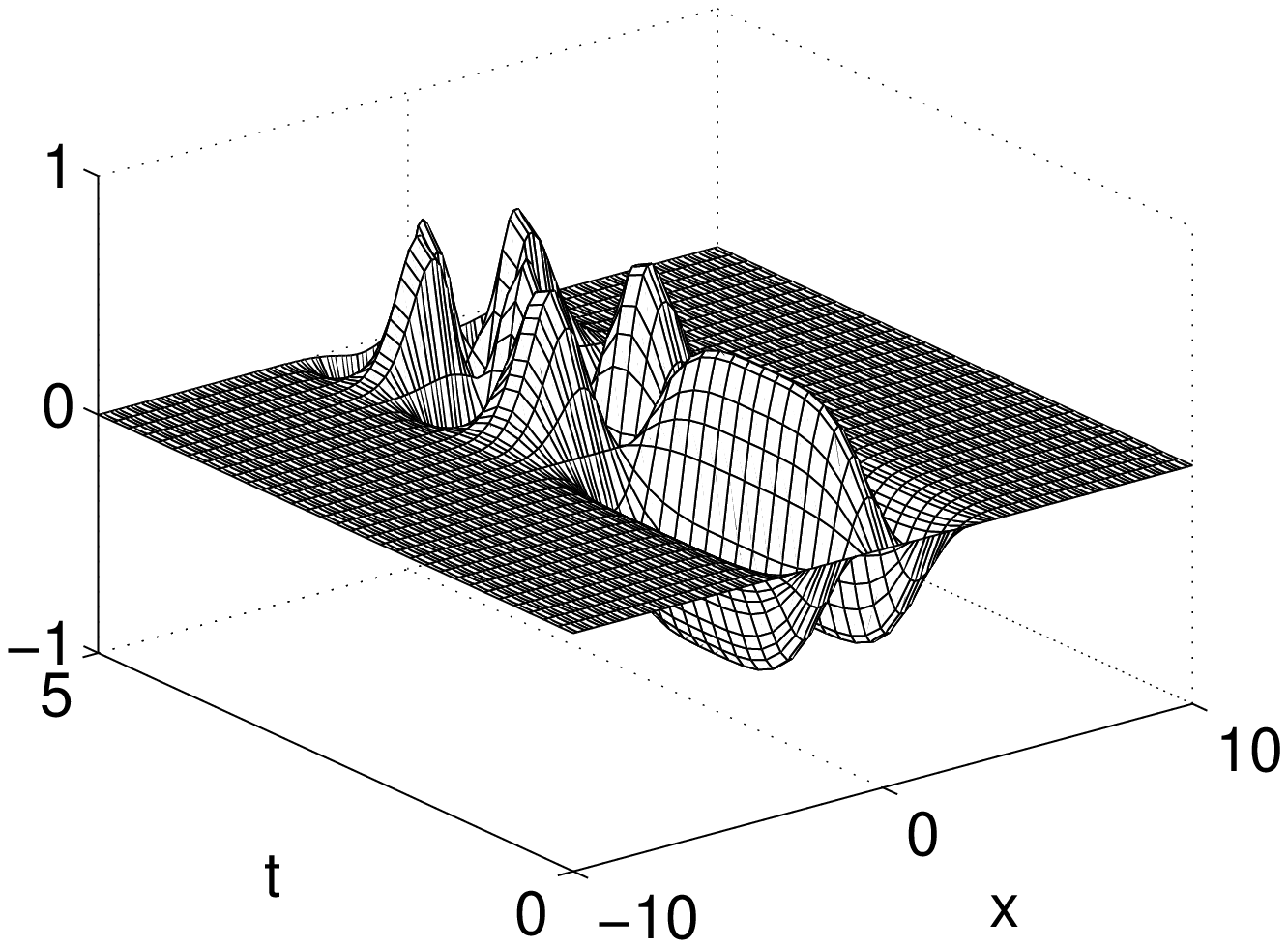}\\
	(a)\\
	\includegraphics[width=6cm]{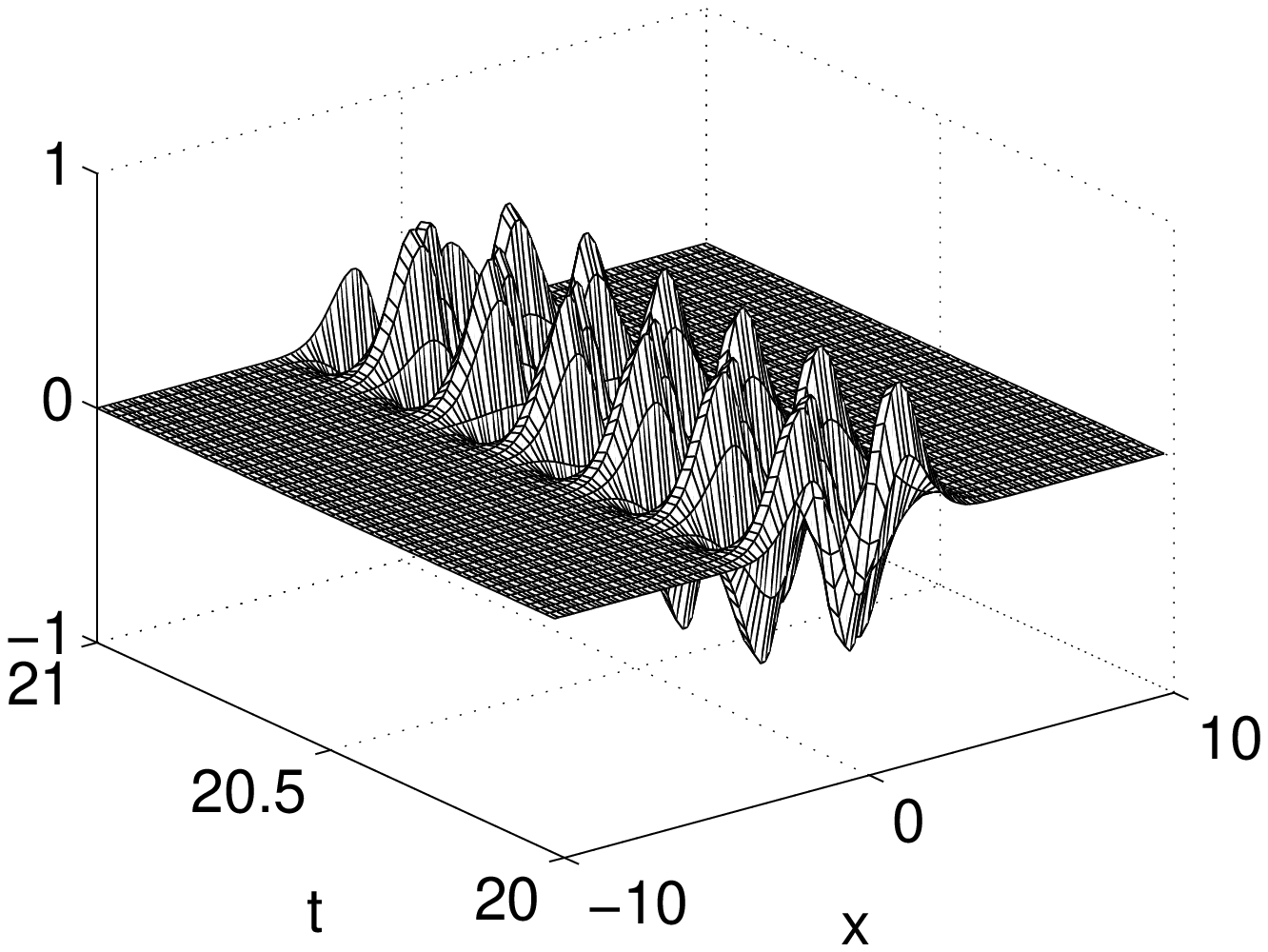}
	\includegraphics[width=6cm]{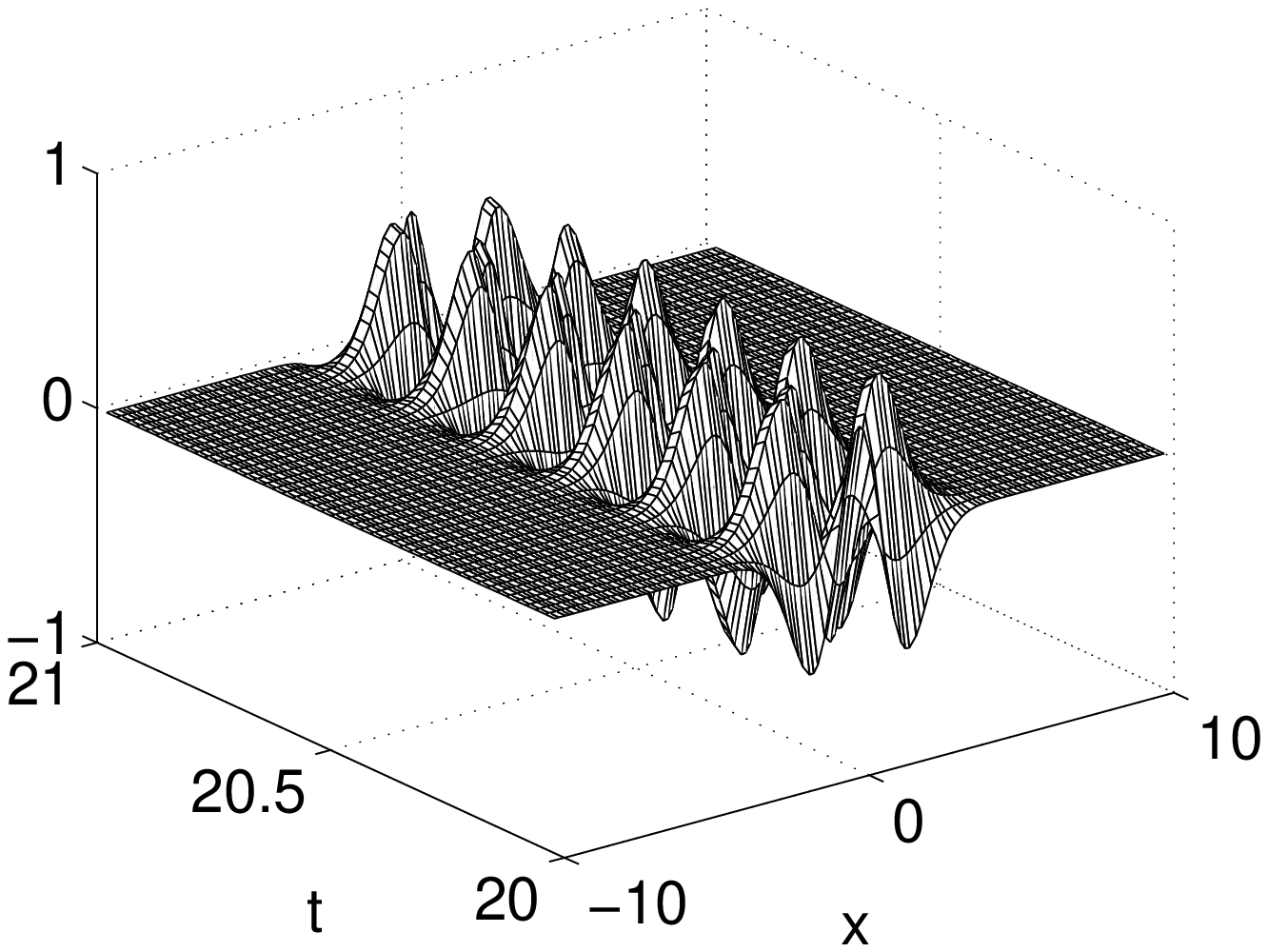}\\
	(b)
				\caption{The real and imaginary part of the CP approximation obtained with $\Delta_x=0.5, \Delta_t=0.25, N=12$ for test problem 1 with $n=2$ in two different $t$-ranges: (a) $t \in[0,5]$,  and (b) the highly-oscillatory region $t\in[20,21]$. }\label{fig:problem1}
\end{figure}
As a first test problem, we consider an equation with a known analytical solution. The potential is given by $V(x,t)=x^2/2-2t$ and $\mu=1$. For the initial wave function we choose an eigenfunction of the harmonic oscillator $\psi(x,0)=(\sqrt{2^n n!\sqrt{\pi}})^{-1/2}\exp(-x^2/2)H_n(x)$ where $H_n$ is the Hermite polynomial of degree $n$. The analytic solution is $\psi(x,t)=\psi(x,0)\exp(-\imath(n+\frac{1}{2})t+\imath t^2)$. The $x$-range we used was $[-10,10]$. We applied our CP-based procedure to approximate $\psi(x,t)$ for $t \in [0,T]$. 
Figure \ref{fig:problem1} gives an idea about the shape of $\psi(x,t)$ when $n=2$. A CP approximation, obtained with the procedure presented in this paper, is shown.
Table \ref{table1} shows some results for different sets of parameter values: $n$ is the number of nodes of the initial wave function, $N$ is the number of terms in the expansion \eqref{eq3}, $\Delta_x$ in the size of the equidistant steps used in the $x$-direction (note that 2 extra inner Lobatto nodes are used in each such step), while $\Delta_t$ is the size of the equidistant time steps. $K$ denotes the number of rectangular sectors. This means that the propagation over each sector requires $T/(K\Delta_t)$ modified Neumann steps. The conservation of norm and the error in the approximations are measured by the following two quantities which are evaluated in $t=T$
\be err_N=\int_{x_{min}}^{x_{max}}\psi^{[K]}(x,T){\psi^{*}}^{[K]}(x,T){\rm d}x-\int_{x_{min}}^{x_{max}}\psi_{exact}(x,T)\psi_{exact}^{*}(x,T) {\rm d}x,\label{errn}\ee
and 
\be err_A=\max_{x\in {\rm mesh}_x}|\psi(x,T)-\psi_{exact}(x,T)|,\ee
where \eqref{errn} is evaluated by applying a classical Lobatto-quadrature rule with the $x$-mesh forming the quadrature nodes. The results in the table illustrate that high accuracies are reached with small $N$, large sector widths and large step sizes. 
For a comparison in efficiency, we included some results obtained with the one-step implicit Crank-Nicolson scheme (CN), which is one of the widespread numerical schemes for solving time-dependent Schr\"odinger problems. The CN scheme requires the solution of a system of algebraic equations of size $(x_{max}-x_{min})/\Delta_x\times (x_{max}-x_{min})/\Delta_x $ at each time-step. It is clear that the CN scheme cannot be used to obtain high accuracies, similar to the ones obtained with the CP-based scheme, within a reasonable time.
All computations were performed in Matlab on a standard desktop computer.

\begin{table}
	\centering
		\begin{tabular}{c|ccc|ccc|cc|c}
		&\multicolumn{3}{c}{$N=10, K=10, T=20$}&\multicolumn{3}{c}{$N=12, K=5, T=20$}&\multicolumn{2}{c}{CN, $T=2$}&CN, $T=5$\\
		\hline
			$\Delta_x$&$1$&$0.5$&$0.25$&$1$&$0.5$&$0.25$&$0.02$&$0.01$&$0.02$\\
			$\Delta_t$&$1$&$1$&$1$&$1$&$1$&$1$&$0.02$&$0.01$&$0.02$\\
			\hline$n=2$&7e-5&1e-10&8e-15&6e-5&6e-11&4e-14&4e-4&1e-4&5e-3\\
			$n=4$&1e-3&3e-10&2e-14&7e-4&1e-10&1e-14&2e-3&5e-4&1e-3\\
			$n=6$&1e-3&5e-10&2e-14&9e-4&2e-10&6e-14&4e-3&9e-4&3e-3\\
			\hline
			$n=2$&1e-4&5e-11&7e-14&5e-4&6e-11&5e-14&3e-4&8e-5&8e-3\\
			$n=4$&5e-4&1e-10&6e-14&5e-4&2e-10&7e-14&1&1&1\\
			$n=6$&8e-4&1e-10&3e-14&6e-4&1e-10&5e-14&8e-1&8e-1&6e-1\\
			\hline
			$time$&6&11&18&4&8&14&174&2450&410
		\end{tabular}
		\caption{Some results for test problem 1. The upper 3 lines of results show $|err_N|$, while the lower ones contain $err_A$. $time$ is the cputime in seconds needed to approximate $\psi(x,t)$ for one $n$ value. The last columns show results obtained by applying the Crank-Nicolson method (CN).}\label{table1}
\end{table}

As a second test problem we used the 1D harmonic oscillator with an explicitly time-dependent frequency \cite{Puzynin}. The Hamiltonian of this system has the form ${\hat H}(x,t)=-\frac{1}{2}\frac{\partial^2}{\partial x^2}+\frac{\omega^2(t)x^2}{2}$ and the frequency $\omega(t)$ was chosen as $\omega^2(t)=4-3 e^{-t}$. The value of the oscillator frequency is quickly doubled over time. The coherent state
\[\psi_0(x)=\left(\frac{1}{\pi}\right)^{1/4} e^{-\frac{1}{2}x^2}\] was chosen as the initial state of the oscillator. Again we choose the segment on $x$ sufficiently large to suppress the influence of the inaccuracy of the boundary conditions on the approximate solution: $x \in [-10,10]$. The error is evaluated in $T=12$. Some results are listed in Tables \ref{table2} and \ref{table2b}.

\begin{table}
	\centering
		\begin{tabular}{c|ccc|cccc|cccc}
		&\multicolumn{3}{c}{$N=10, K=20$}&\multicolumn{4}{c}{$N=15, K=20$}&\multicolumn{4}{c}{$N=20, K=20$}\\
		\hline
			$\Delta_x$&$1$&$0.5$&0.25&$1$&0.5&0.25&0.2&0.5&0.25&0.2\\
			$\Delta_t$&$0.6$&$0.3$&0.05&$0.6$&0.3&0.05&0.02&0.3&0.05&0.02\\
			\hline
			$|err_N|$&2e-3&5e-5&3e-9&2e-3&5e-5&6e-9&6e-11&6e-6&6e-9&6e-11\\
			$err_A$&3e-3&8e-5&1e-5&3e-3&8e-5&3e-7&3e-7&8e-5&6e-8&2e-9\\
			\hline
			$time$&10&24&62&17&36&104&258&55&184&390
		\end{tabular}
		\caption{Error estimations for the CP-approach applied on test problem 2. $time$ is the CPU time in seconds needed to approximate $\psi(x,t)$ over $[-10,10]\times[0,12]$ using stepsizes $\Delta_x$ and $\Delta_t$.}\label{table2}
\end{table}

\begin{table}
	\centering
		\begin{tabular}{c|ccc}
		$N$&$K=1$&$K=10$&$K=20$\\
		\hline
			5&1e-2&3e-3&2e-3\\
			10&3e-4&3e-5&9e-6\\
			15&2e-5&3e-6&3e-7\\
			20&7e-7&4e-8&2e-9\\
			30&2e-9&3e-11&1e-12\\
			\hline
		\end{tabular}
		\caption{Error estimations $err_A$ for test problem 2. Different values of $K$ and $N$ were used to approximate $\psi(x,t)$ over $[-10,10]\times[0,12]$ using stepsizes $\Delta_x=0.2$ and $\Delta_t=0.02$.}\label{table2b}
\end{table}

As third test problem, we consider the Walker and Preston model of a diatomic molecule in a strong laser field \cite{Walker}. This system is described by the one-dimensional Schr\"odinger equation
\[i\frac{\partial \psi(x,t)}{dt}=\left(-\frac{1}{2\mu}\frac{\partial^2}{\partial x^2}+V(x)+f(t)x\right)\psi(x,t),\]
with $\psi(x,0)=\psi_0(x)$. Here $V(x)=D(1-e^{-\alpha x})^2$ is the Morse potential and $f(t)x=A\cos(\omega t)x$ accounts for the laser field. This problem is used as a test bench for the time propagation methods presented in e.g. \cite{Gray} and \cite{JM}. As in \cite{Gray} we take values for the parameters (in atomic units: a.u.) corresponding to the HF molecule: $\mu=1745$ a.u., $D=0.2251$ a.u., $\alpha=1.1741$ a.u., $A=0.011025$ a.u. and laser frequency $\omega=0.01787$. We assume the system is defined in the interval $x\in [-1,4.32]$. As initial conditions we take the ground state of the Morse potential $\psi_0(x)=\sigma \exp(-(\rho-1/2)\alpha x)\exp(-\rho e^{-\alpha x})$, with $\rho=2D/\omega_0$, $\omega_0=\alpha\sqrt{2D/\mu}$ and $\sigma$ the normalizing constant, i.e. $\sigma=0.2411580885\times 10^{-10}$.

A natural time unit for the problem is $\tau=2\pi/\omega$. In order to illustrate that the approach can also be used over longer time ranges, the integration was carried out over the time range $t\in[0,100\tau]$. The exact solution is accurately approximated using a sufficiently small time step. As in \cite{Gray,JM}, we choose an equidistant grid for the spatial coordinate $x$ with 64 steps. The results for different time-steps $\Delta_t$ are shown in Table \ref{table3}. These results were obtained with $N=15, K=20$.
When taking a fixed set of basis functions, e.g.\ $N$ eigenfunctions of the time-independent harmonic oscillator, over the entire $t$-range, similar accuracies can only be reached by using values for $N>100$.

\begin{table}
	\centering
		\begin{tabular}{c|cccccc}
		\hline
			$\Delta_t$&$\tau/10$&$\tau/20$&$\tau/50$&$\tau/100$&$\tau/200$&$\tau/400$\\
			\hline
			$|err_N|$&3e-3&1e-4&1e-6&2e-9&5e-10&5e-10\\
			$err_A$&1e-2&7e-4&2e-5&1e-6&8e-8&5e-9\\
			\hline
			$time$&114&183&412&863&2041&5652
		\end{tabular}
		\caption{Error estimations for the CP-approach applied on test problem 3. We used $N=15, K=20$.}\label{table3}
\end{table}

\section{Conclusion}
To accurately simulate the outcome of quantum dynamical experiments with time-dependent Hamiltonians, accurate and efficient numerical techniques should be used. In this paper, we showed how the succesful CP technique can be used in the spatial discretization of a time-dependent Schr\"odinger equation following an approach suggested by Ixaru in \cite{Ixaru2010}, but also for efficient and accurate time stepping in the resulting ODE system. The effectiveness of the approach was illustrated by some numerical examples.

\providecommand{\Gr}{Gr} \providecommand{\ea}{et al}

\end{document}